\documentclass[letterpaper, 10 pt, conference]{ieeeconf}  

\IEEEoverridecommandlockouts                              

\usepackage{amsmath,amssymb,amsfonts}

\newcommand{\R}{\mathbb{R}}

\newtheorem{theorem}{Theorem}[section]

\usepackage{graphicx,color}

\include{Xlatin1}

\newcommand{\dis}{\displaystyle}

\newcommand{\be}{\begin{equation}}
\newcommand{\ee}{\end{equation}}
\newcommand{\ba}{\begin{array}}
\newcommand{\ea}{\end{array}}

\newcommand{\hsm}{\hspace{-.1cm}}
\newcommand{\vsm}{\vspace{-.1cm}}

\title{\LARGE \bf
Crowd motion paradigm modeled by a bilevel sweeping control problem}

\author{Tan H. Cao$^{1}$, Nathalie T. Khalil$^{2}$, Boris S. Mordukhovich$^{3}$, Dao Nguyen$^{3}$, and Fernando Lobo Pereira$^{2}$
\thanks{$^{1}$T.H. Cao is with Department of Applied Mathematics and Statistics,	SUNY (State University of New York) Korea
	{\tt\small tan.cao@stonybrook.edu}}%
\thanks{$^{2}$N.T. Khalil and F.L. Pereira are with SYSTEC, Faculty of Electrical Engineering, Porto University, and with the Institute for Systems and Robotics, 4200-465 Porto, Portugal
        {\tt\small khalil.t.nathalie@gmail.com, flp@fe.up.pt}}%
    \thanks{$^{3}$B.S. Mordukhovich and D. Nguyen are with Department of Mathematics, Wayne State University, USA
    	{\tt\small boris@math.wayne.edu, dao.nguyen2@wayne.edu}}%
}
\begin{document}

\maketitle
\thispagestyle{empty}
\pagestyle{empty}
\begin{abstract}
This article concerns an optimal crowd motion control problem in which the crowd features a structure given by its organization into $N$ groups (participants) each one spatially confined in a set. The overall optimal control problem consists in driving the ensemble of sets as close as possible to a given point (the ``exit'') while the population in each set minimizes its control effort subject to its sweeping dynamics with a controlled state dependent velocity drift. In order to capture the conflict between the goal of the overall population and those of the various groups, the problem is cast as a bilevel optimization framework. A key challenge of this problem consists in bringing together two quite different paradigms: bilevel programming and sweeping dynamics with a controlled drift. Necessary conditions of optimality in the form of a Maximum Principle of Pontryagin in the Gamkrelidze framework are derived. These conditions are then used to solve a simple illustrative example with two participants, emphasizing the interaction between them.
\end{abstract}

\section{Introduction}\label{sec:intro}

{Problems modeled as a bilevel optimization, with dynamics featuring a sweeping process control arise naturally in numerous applications. For instance in managing the motion of structured crowds organized in groups, in operating teams of drones providing complementary services
	in a shared confined space, in nanoferro-electric technologies for functional improvement of mobile electronic devices, among many others. In all those examples, the problem can be represented as a bilevel optimization, with dynamics modeling some structure formed by a set of groups with distinct properties and confined to controlled bounded moving subsets, giving rise to a sweeping process control phenomena.}

The purpose of this article is to present a bilevel sweeping control problem through a model arising in the management of structured crowd motions on the plane. This framework, bringing together bilevel optimization and controlled sweeping processes, was addressed for the first time in \cite{KhalilPereiraCDC2019}. By structured crowd we refer to a {\it population} organized into {\it group of sets}. As an example, let us imagine a population trying to exit a certain space with the shortest possible path. This population is organized into groups of sets, each moving  along a trajectory $y_i$ prescribed by a coordinator in order to reach the exit. The population in each group has to remain inside its own moving set, while minimizing its effort to achieve this. In this context, we formulate a bilevel problem coupled with a sweeping control process: the upper level defines the direction of each group (or set) with the goal of driving the ensemble of sets as close as possible to the exit, while avoiding any overlapping between them, and the lower level problem where each group population has to stay confined to its moving group, via a motion modeled by a sweeping control process, while minimizing its control effort to achieve this.

\begin{figure}[h!]\begin{center}
\includegraphics[width=6.8cm, height=4.5cm]{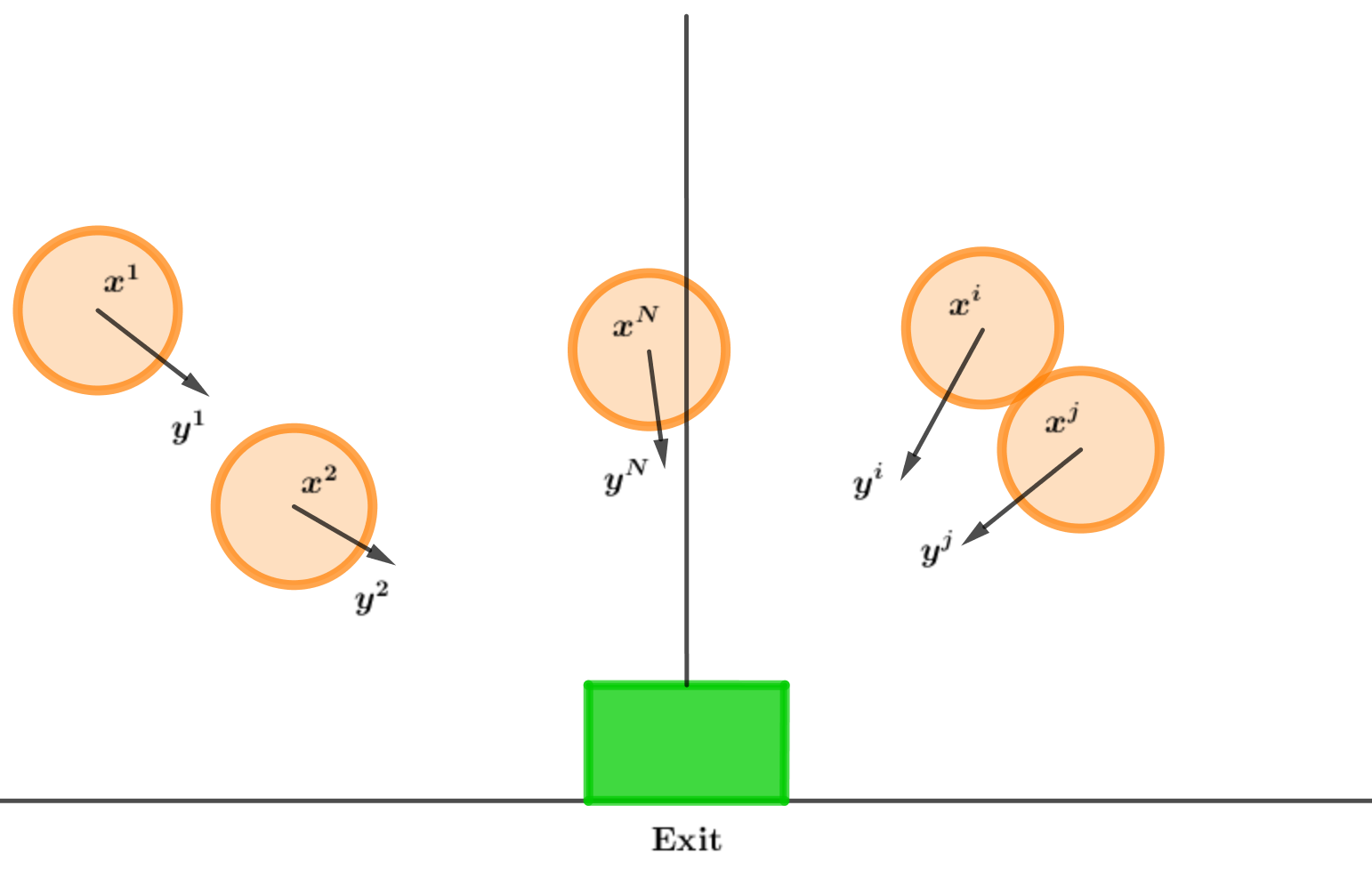}
\vspace{-.3cm}\caption{Crowd motion model for $N$ participants}\label{pic:general case}
\end{center}
\end{figure}

To simplify, we take the groups to be $N$ disks, in the plane, of the same radius $R$. Each disk is subject to a translation vector $y^i\hsm\in\hsm\R^2$ {representing its linear motion direction}. The population in each disk  will be presented by its ``representative'' position, $x^i\hsm\in\hsm \R^2$. The exit set the origin. The problem is illustrated in Figure \hsm\ref{pic:general case}.

\vspace{-.4cm}
\vspace{-.5cm}
To formalize the ideas stated above, for $i=1,\ldots,N$, we denote by $(P_H(x_0,u))$ the upper level problem in $\R^2$:
\begin{eqnarray}(P_H(x_0,u)) &&\hspace{-.4cm} \mbox{Minimize}\;  J_H(y;x_0,u) \nonumber\label{highcost}\\
	\mbox{subject to }&& \hspace{-.4cm} \dot y^i(t)= v^i(t)\quad [0,T]\mbox{-a.e.}\label{highdynam}\nonumber, \quad  y^i(0)=y_0^i \;\label{highendp}\nonumber\\
	&&  \hspace{-.4cm} v^i \in \mathcal{V}^i\hsm :=\hsm  \{v^i\hsm\in\hsm L^2([0, T]; \R^2): v^i(t)\hsm \in\hsm V^i\}\nonumber\\
	&& \hspace{-.4cm}  \|y^i(t)-y^j(t)\| \ge 2R, \; j\neq  i \; \forall t\label{condition:nonoverlaping}\\
	&& \hspace{-.4cm} \mbox{and } y^i \mbox{ s.t. } \exists \mbox{ solution } { (x_0^i,u^i)} \text{ to } P_L^i(y^i) \vsm\vsm\label{low-level feasibility}
\end{eqnarray}
with $(y,x_0,u)^T\hsm \in\hsm \R^{N(4+m)}$, $a\hsm =\hsm\textit{col}\{a^i\}_{i=1}^N$, \textit{``col''} stands for ``column'', and $\dis J_H(y;x_0,u)\hsm:= \hsm\frac{1}{2} \sum_{i=1}^{N} \|y^i(T)\|^2$. $T$, and $y_0^i$ are fixed, and for $i,j=1,\ldots,N$ with $j\hsm\neq\hsm i$, $\|y_0^i\hsm -\hsm y_0^j\|\hsm \geq\hsm 2R$, $\|\cdot\|$ is the Euclidean norm in $\R^{2}$. $V^i \hsm \subset \hsm \R^{2} $ stands for the control set and is compact. Condition (\ref{condition:nonoverlaping}) represents the non-overlapping imposed on the translated disks $D^i(t) \hsm:=\hsm D\hsm +\hsm y^i(t)$, $ D^j(t) \hsm:=\hsm D\hsm +\hsm y^j(t)$ with $D \hsm =\hsm\{x\hsm \in\hsm \R^2\hsm :\|x\|\hsm\leq\hsm R\}$, i.e., the distance between the centers of two disks have to be not less than $2R$.

The problem $(P_L^i(y^i))$, alluded to in (\ref{low-level feasibility}), is defined by the following parametric lower level problem with dynamics involving a sweeping process with a controlled drift.
\begin{eqnarray}
(P_L^i(y^i)) &&\hspace{-.4cm}\mbox{Minimize}\; J_L^i(x_0^i,u^i;y^i) \label{cost-i}\nonumber\\
\mbox{subject to }&& \hspace{-.4cm}\dot x^i(t)\hsm\in\hsm f^i(x^i(t),u^i(t))\hsm-\hsm N_{D+y^i(t)}^{M^i} (x^i(t)) \label{sweep-i}\label{eq:dynamic lower original} \\
&&\nonumber\hspace{4.5cm}  [0,T]\mbox{-a.e.}\vsm\vsm\vsm \\
&&\hspace{-.4cm} x^i(0)=x_0^i \in D+y_0^i\nonumber\\
&&\hspace{-.4cm} u^i(t)\in\mathcal{U}^i\hsm :=\hsm \{ u^i\hsm\in\hsm L^\infty([0,T];\R^m) \hsm: \hsm u^i(t)\hsm\in\hsm U^i\}\nonumber\\
&&\hspace{-.4cm} x^i(t)\in D+y^i(t)\quad \forall t\in [0,T],\vsm\vsm\label{concons-i}
\end{eqnarray}
where, for a given process $y$ solving $(P_H(x_0,u))$, $\dis J_L^i(x_0^i,u^i;y^i) := \int_{0}^{T}\hsm \|u^i(t)\|^2 dt$ (control effort of each disk population), $f^i\hsm :\hsm\R^{2}\hsm\times\hsm \R^{m} \to\R^{2}$, $U^i\hsm\subset\hsm \R^{m}$ is compact, the truncated cone $ N_{A^i}^{M^i}(z)\hsm :=\hsm N_{A^i}(z) \cap  M^i B_1(0)$, being $N_{A^i}(z)$ the Mordukhovich (limiting) normal cone to the closed set $A^i$ at point $z$ in the sense of \cite{mordukhovich2006book}, $ B_1(0)$ the closed unit ball in $\R^2$ centered at the origin, and $M^i\hsm >\hsm 0$ a given constant.

Let us define some concepts. Let $(P_L)\hsm :=\hsm \textit{col}\{(P_L^i)\}_{i=1}^N$. For some parameter $y$, a pair $(x(\cdot), u(\cdot))$ is a {\it feasible} (or admissible) control process to $(P_L)$ if $u(\cdot)$ is feasible control to $(P_L)$, and $x(\cdot)$ is an arc satisfying the differential inclusion (\ref{sweep-i}), the initial condition, together with (\ref{concons-i}). An {\it optimal solution} to $(P_L)$ is a feasible pair of $(P_L)$ minimizing the value of the cost functional $J_L(x_0,u;y)$ over all admissible pairs of $(P_L(y))$. A feasible triple of the dynamic control problem $(P_H)$ is the set of feasible process $y$, and an optimal pair $(x_0,u)$ to $(P_L(y))$. The feasible triple $(y;x_0,u)$ is optimal to $(P_H)$ if $(y;x_0,u)$ minimizes the value of $J_H(y;x_0,u)$ among all admissible strategies of $(P_H)$.

It is important to note that the articulation of $(P_H)$, and $(P_L)$ in this article features significant differences with respect to the bilevel control sweeping process problem studied in \cite{KhalilPereiraCDC2019}, and in \cite{KhalilPereiraIFAC2020}. Here, the upper level problem $(P_H)$ acts on the dynamics of each lower level sweeping control problems $(P_L^i)$, via the arcs $y^i(\cdot)$ specifying the motion of the respective disks $D$. While the upper level problem minimizes the sum of distances of all the disk centers to the exit at the final time, and forbids their overlap, the lower level problem minimizes the control effort spent by any disk population in order to remain within the constraint set $\dis \Pi_{1}^ND+y(t)$. Here $y$ is a feasible arc to $(P_H)$, and a parameter to $(P_L)$, and $\Pi_{1}^N D $ denotes the Cartesian product of $N$ disks $D$. We establish necessary conditions of optimality for this problem, and use them to solve an example with $N\hsm =\hsm 2$.

In addition to the fact that a time-optimal problem is considered at the upper level instead of the minimal-path studied in the current article, in \cite{KhalilPereiraCDC2019}, the $N$ disks are confined in a larger constraint set, and $y^i$, $i=1,\ldots,N$ can take values on the boundary of this extra constraint set, giving rise to a sweeping process at the upper level. Only well-posedness and the existence of solutions to the problem are studied in \cite{KhalilPereiraCDC2019}. However, in \cite{KhalilPereiraIFAC2020}, a simpler instance of the problem n \cite{KhalilPereiraCDC2019} is studied as only one disk (i.e. $N=1$) intervenes, and no sweeping process appears at the upper level problem. Necessary conditions of optimality are established for this problem.

The problem in \cite{CaoMordukhovich2018DCDS} is of a different nature as it considers only a single-level (and not a bilevel) control sweeping problem. Another main difference is the nature of the sweeping process arising in the dynamics. Indeed, while in our paper the sweeping process appears, at the lower level problem, to force the various populations, with (average) motion velocity $\dot x^i(t)$ to stay confined to their moving set $D+y^i(t)$, in \cite{CaoMordukhovich2018DCDS} the sweeping process intervenes to adjust velocities when a contact occurs between the disks, supposed to remain at a minimum distance from each other, avoiding overlapping scenarios. On the other hand, the technique used to deal with the sweeping process in \cite{CaoMordukhovich2018DCDS} relies on the discrete approximation, while in our paper, we adopt another technique by approximating the sweeping term by a sequence of Lipschitz functions in the state variable. More details can be found in \cite{KhalilPereiraIFAC2020}, and in Section \ref{sec:proof} of this paper.

\vspace{.1cm}
{\bf Notation.}  We denote by $N_A(x)$ the Mordukhovich (limiting) normal cone to a closed set $A$ at the point $x\in A$, and by $\partial \varphi$, and $\partial^C \varphi$, respectively, the Mordukhovich, and the Clarke subdifferentials of $\varphi$. If $\varphi$ is locally Lipschitz, then $\partial^C \varphi= \textrm{co}\ \partial \varphi$, where ``$\textrm{co} A$'' denotes the closure of the convex hull of the set $A$. We refer the reader to \cite{mordukhovich2006book,Clarke1990,vinter2010book} for a full overview on nonsmooth analysis. $AC([0,T];\R^{2N})$ stands for the space of absolutely continuous functions on $[0,T]$ with values in $\R^{2N}$, $BV([0,T];\R^+)$ for the space of functions of bounded variations on $[0,T]$ with nonnegative values, $\|\cdot\|_{TV}$ the total variation, and superscript $T$ the transpose. {Given $ a\hsm\in\hsm\R^k $, and $ b\hsm\in\hsm\R^{km} $, we have $ a\diamond b \hsm =\hsm {\it col}\left(a^i(b^{(i-1)m+1},\cdots, b^{im})^T: i\hsm = \hsm 1,\ldots,k\right)$.}

The article is organized as follow: in Section \ref{sec:assumptions}, we give the required assumptions to be imposed on the data. In Section \ref{sec:NCO}, we establish the necessary optimality condition of the problem setting $(P_H)$-$(P_L)$. A proof outline with the key ideas is given in Section \ref{sec:proof}. An example considering the case of only two disks is studied in Section \ref{sec:examples}. We finish with a conclusion and some remarks for future avenues.

\section{Assumptions}\label{sec:assumptions}
Before stating the necessary optimality conditions, we present the assumptions to be imposed on the data of the problem. These are as follows, for all $i=1,\ldots,N$:
\begin{itemize}
\item[H1] $f^i(x^i,\cdot)$ is Borel measurable $\forall\,x^i\hsm\in\hsm \R^2$, $f^i(\cdot,u^i)$ is Lipschitz continuous $\forall\,u^i\hsm\in\hsm U^i$, and $f^i(\cdot,\cdot)$ is bounded $\forall (x^i,u^i)\hsm \in\hsm \R^2\hsm \times\hsm U^i$.
\item[H2] $f^i(x^i,U^i)\subset\R^2$ is a closed and convex set for each $x^i$.
\item[H3] The control sets $U^i$, and $V^i$ are compact and convex.
\item[H4] There exists $\beta \hsm>\hsm 0$ s.t. $\beta B_1(0)\hsm \subset\hsm f^i(x^i,U^i),\, \forall\, x^i\hsm\in\hsm \R^2$.
\item[H5]  The constant $M^i$ specifying the truncation of the normal cone satisfies $\overline M^i\hsm >\hsm M^i\hsm >\hsm \overline m^i$
where, $ \forall \zeta^i\hsm \in\hsm N_{D+y^i(t)}(x^i(t)), \forall t\hsm\in\hsm[0,T]$ with $ x^i(t)\hsm\in\hsm \text{bd} (D+ y^i(t))$ {(``bd'' is the boundary)},
\vsm\vsm
\begin{equation}
\hspace{-.3cm}\overline M^i \hsm=\hsm \min_{\|\zeta^i\|=1}\left\{\max_{u\in U^i}\{\langle \zeta^i, f^i(x^i(t),u)\rangle\} -\hsm\min_{v\in V^i}\{\langle\zeta^i, v\rangle\}\right\},\vsm\vsm \nonumber
\end{equation}
\begin{equation}
\hspace{-.3cm}\overline m^i\hsm =\hsm \max_{\|\zeta^i\|=1}\left\{\min_{u\in U^i}\{\langle \zeta^i, f^i(x^i(t),u)\rangle\}- \max_{v\in V^i}\{\langle\zeta^i, v\rangle\}\right\}.\nonumber
\end{equation}
\end{itemize}
We also require additional assumptions playing a critical role in the articulation of  $(P_H)$, and $(P_L)$, and in the derivation of the necessary conditions of optimality. Similar assumptions has been cpnsidered in \cite{Ye1997,BenitaDempeMehltiz2016}, albeit for in a different context:
\begin{itemize}
\item[H6] For a given $y$ feasible to $(P_H)$, the solution set of $(P_L)$ is not empty and every process $(x,u)$ solving $(P_L)$ lies in the interior of its solution set.
\item[H7] The articulation of $(P_H)$-$(P_L)$ is partially calm, i.e. , $\exists\,\rho^i \ge 0$ such that for any feasible $(y,x,v,u)$,
\vsm
\begin{align}
&\nonumber J_H(y;x_0,u) - J_H(y^{*};x^{*}_0,u^{*})\vsm\\ &\hspace{1cm}+ \sum_{i=1}^{N}\rho^i \hsm\left(\int_0^T\hsm\hsm \|u^i(s)\|^2ds-\varphi^i(v^i)\hsm\right) \geq 0,\label{def:partial calmness}
\end{align}  where $(y^{*},x^{*}_0,u^{*})$ is an optimal solution to $(P_H)$-$(P_L)$, and $\varphi^i(\cdot)$ is the value function of $(P_L^i)$ defined by
\vsm
\begin{align}\hsm\varphi^i(v^i) \hsm  =\hsm  \nonumber\min \Bigg\{\hsm J_L^i&(x_0^i,u^i;y^i)  \hsm :\hsm (x_0^i,u^i) \mbox{ feasible for } (P_L^i)\vspace{-.5cm} \\ & \hspace{.2cm}\mbox{ and } y^i(t) = y_0^i+\hsm \int_0^t \hsm v^i(s)ds\hsm \Bigg\}, \vsm\vsm\vsm\label{def:value function-low} \end{align}
being $v^i$ the $i^{th}$ component of $v$, an admissible to $(P_H)$.
\end{itemize}

\section{Necessary optimality conditions}\label{sec:NCO}

Before stating the necessary optimality conditions of $(P_H)$-$(P_L)$, we shall first define the following\vsm
\begin{eqnarray*} && \hspace{-.7cm} H_H(y,x,v,u,q_H,q_L,\nu_H,\nu_L, \alpha) \\ && = \sum_{i=1}^N H_H^i(y^i,x^i,v^i,u^i,q_H^{i},q_L^i,\bar \nu_H^{i},\nu_L^i, \alpha^i) \end{eqnarray*}
where $y$, $ x$, $v$, $u$, $q_H$, $q_L $ take values in $\R^{2N}$, $\nu_L$ are in $\R^{N}_+$, $\nu_H \hsm \in\hsm \R^{N(N-1)}_+\hsm$, $\bar \nu_H^i\hsm \in\hsm \R^{N-1}_+$ is a vector with components $\nu_H^{ij}$ satisfying $j\hsm\neq\hsm i$, and $\nu_H^{ij}\hsm=\hsm\nu_H^{ji}$ (symmetric), and for a fixed $i \hsm\in \hsm \{1,\ldots,N\}$, we have\vsm
\begin{eqnarray*}
&&\hspace{-.7cm} H_H^i(y^i,x^i,v^i,u^i,q_H^{i}, q_L^i,\bar \nu_H^{{i}},\nu_L^i, \alpha^i) \\ &&:=\langle q_L^i- \nu_L^i(x^i- y^i),f^i(x^i,u^i)\rangle + \nu_L^i \langle x^i- y^i,v^i\rangle\hsm \\
&&\hspace{.7cm}+\sigma^i(y^i,x^i,q_L^i,\nu_L^i) - \alpha^i\|u^i\|^2 \vsm \\
&& \;\;+ \left \langle q_H^i\hsm +\hsm \sum_{j<i}^{N}\hsm \nu_H^{ij}\frac{y^i-y^j}{\|y^i-y^j\|}+\hsm \sum_{j>i}^{N}\hsm \nu_H^{ij}\frac{y^i-y^j}{\|y^i-y^j\|} ,v^i\hsm\right\rangle,
\end{eqnarray*}
where $\dis \sigma^i(y^i,x^i,q_L^i,\nu_L^i) = \hsm\hsm\hsm\sup_{\xi\in-N^{M^i}_{D^i}(x^i-y^i)} \hsm\hsm\hsm \{\langle q_L^i\hsm -\hsm\nu_L^i(x^i\hsm - \hsm y^i),\xi\rangle\}$.

We note that the multiplier $\bar \nu_H^i= \{\nu_H^{ij}\}_{i\ne j}$ appears in $H_H^i$ to reflect the activity of the constraint  $\|y^i\hsm \hsm -y^j\| \hsm\geq\hsm 2R$ (i.e., (\ref{condition:nonoverlaping})). It is non-increasing whenever $\|y^i\hsm -\hsm y^j\|\hsm =\hsm 2R$ (i.e., the disks $D\hsm +\hsm y^i(t)$ and $D\hsm +\hsm y^j(t)$ are in contact), and constant otherwise.

The form of the Hamiltonian $H_H$ stated above is different from the usual Pontryagin-Hamilton function used in the Dubovitskii-Milyutin form \cite{Dubovitskii 1965}. Our Hamiltonian is the one used in establishing the necessary conditions of optimality in the Gamkrelidze's
form \cite{Gamkrelidze1959,Gamkrelidze1960}. These results, and further developments were recently incorporated in the modern optimal control literature in \cite{Arutyunov_Karamzin_Pereira_2011,Arutyunov-Karamzin-Pereira2017,Karamzin_Pereira_2019}. This form differs from the Dubovitskii-Milyutin one in the way the measure multiplier associated with the state constraints enters in the Pontryagin-Hamiltonian function. It might entail some loss of generality due to the extra smoothness required on the function specifying the state constraints, but opens significant new computational perspectives for indirect methods based on the Maximum Principle of Pontryagin due to the regularity of its measure multiplier.

\begin{theorem}\label{theorem:main thm} Let H1-H7 hold and $(y^*,x^*,u^*)$ a solution to $(P_H)$-$(P_L)$. Then, there exists a set of  multipliers $(q_H,q_L,\nu_H,\nu_L,\lambda, \alpha)$ with $q_H$, and $q_L$ in $AC([0,T];\R^{2N})$, $\nu_H\hsm\in \hsm BV([0,T];\R^{N(N-1)}_+)$, and $\nu_L \hsm\in\hsm BV([0,T];\R^{N}_+)$ being both non-increasing, and $\nu_H^{ij}$, and $\nu_L$ constants on $\{t\hsm\in\hsm[0,T]\hsm:\hsm \|y^i-y^j\|\hsm >\hsm 2R, j\neq i\}$, and $\{t\hsm\in\hsm[0,T]\hsm:\hsm \|y-x\|\hsm <\hsm R\}$, respectively, and $\lambda \in [0,1]$, $\alpha\hsm\in\hsm\R^N$, with $\alpha^i=\lambda\rho^i, i=1,\ldots,N$ (being $\rho^i$ the modulus in (H7)), with:
\begin{itemize}
\item[1.] Nontriviality.
$\dis \hsm\|(q_H,q_L)\|_{L^\infty}\hsm+\hsm \|(\nu_H,\nu_L)\|_{TV} \hsm+\hsm\lambda \hsm +\hsm |\alpha| \hsm\neq\hsm  0$\vspace{.1cm}
\item[2.] Adjoint equations.
\begin{eqnarray*}
-\dot q_L(t) &\hsm\hsm\in&\hsm\hsm \partial_x H_H (y^*,x^*,v^*,u^*,q_H,q_L,\nu_H,\nu_L, \alpha)  \nonumber\\
&\hsm\hsm=&  \hspace{-.1cm} \hsm\hsm\partial_x\langle q_L(t)\hsm-\hsm {\nu_L(t)\diamond (x^*(t)\hsm-\hsm y^*(t))},f^*(t) \rangle \nonumber \\
&& \hspace{-1.3cm} + \nu_L(t) v^*(t) \hsm +\hsm \partial_x \sigma(y^*(t),x^*(t),q_L(t),\nu_L(t))  \; \mbox{a.e.}\nonumber\vspace{.2cm} \\ 			
-\dot q_H(t) &\hsm\hsm\in\hsm\hsm &\partial_y H_H (y^*,x^*,v^*,u^*,q_H,q_L,\nu_H,\nu_L, \alpha) \nonumber  \\
& \hsm\hsm =\hsm\hsm& -{\nu_L(t)\diamond v^*(t) + \nu_L(t)\diamond f^*(t)}\nonumber\\
&& \hspace{.4cm}+ \partial_y \sigma(y^{*}(t),x^{*}(t),q_L(t),\nu_L(t))\nonumber\\
&&  \hspace{0.4cm}+ \textit{col}\bigg(\sum_{j>1}   \nu_H^{{1j}}(t) d^{1j}(t)v^{1*}(t);\vsm \\
&&\hspace{1cm}\ldots;\nonumber  \sum_{j<i}  \nu_H^{{ij}}(t) d^{ij}(t)v^{i*}(t)\vsm \\
&& \hspace{1.2cm} + \sum_{j>i} \nu_H^{{ij}}(t) d^{ij}(t)v^{i*}(t); \nonumber	 \\
&& \hspace{.6cm} \ldots; \sum_{j<N} \hsm \nu_H^{{Nj}}(t) d^{Nj}(t)v^{N*}(t) \bigg) \mbox{ a.e.}\vsm\vsm\nonumber			
\end{eqnarray*}
being $\dis d^{ij}\hsm:=\hsm \|y^{i*}\hsm -\hsm y^{j*}\|^{-1} I\hsm -\hsm \frac{(y^{i*}\hsm -\hsm y^{j*})(y^{i*}\hsm -\hsm y^{j*})^T}{\|y^{i*}-y^{j*}\|^{3}}$, $I$ the unit matrix, $f^*(t)\hsm:=\hsm f(x^*(t),u^*(t))$, and $\sigma \hsm=\hsm\sum_{i=1}^N \sigma^i.$\vspace{.1cm}
\item[3.] Boundary conditions.
\vspace{-.2cm}
\begin{eqnarray*}
q_H(0) &\hsm\hsm\in\hsm\hsm& \R^{2N} \\q_L(0)&\hsm\hsm\in\hsm\hsm&  N_{D+y_0}(x^*(0)) + {\nu_L(0)\diamond (x^*(0)-y_0)}\\
q_H(T)&\hsm\hsm =\hsm\hsm& -\lambda y^*(T) -{\nu_L(T) \diamond  (x^*(T)-y^*(T))}\\
&& \hspace{.2cm} - \textit{col}\Bigg(\sum_{j>1}  \nu_H^{{1j}}(T) \frac{y^{1*}(T)-y^{j*}(T)}{\|y^{1*}(T)-y^{j*}(T)\|} ; \vsm\\
&& \hspace{.6cm}\ldots; \sum_{j<i}   \nu_H^{{ij}}(T) \frac{y^{i*}(T)-y^{j*}(T)}{\|y^{i*}(T)-y^{j*}(T)\|} \\
&& \hspace{.8cm}+ \sum_{j>i}   \nu_H^{{ij}}(T) \frac{y^{i*}(T)-y^{j*}(T)}{\|y^{i*}(T)-y^{j*}(T)\|}; \\
&&\hspace{.3cm} \ldots; \sum_{j<N}  \nu_H^{{Nj}}(T) \frac{y^{N*}(T)-y^{j*}(T)}{\|y^{N*}(T)-y^{j*}(T)\|} \Bigg)\nonumber \\
q_L(T)  &\hsm\hsm=\hsm\hsm& {\nu_L(T)\diamond (x^*(T)-y^*(T))}.\vsm\vsm\vsm
\end{eqnarray*}
\item[4.] Maximum condition on the lower level control.
		
$u^{*}(t)$ maximizes on $U\hsm =\hsm U^1\hsm\times\hsm\ldots\hsm\times\hsm U^N$, $[0,T]$-a.e.
\vspace{-.1cm}
\begin{eqnarray*}(u^1,\ldots,u^N)\to \hspace{-.3cm}  &&\hsm\hsm\hsm - \sum_{i=1}^{N}\alpha^i \|u^i\|^2\vsm\\
&&\hspace{-1.5cm} +\langle q_L(t)\hsm -\hsm  {\nu_L(t)\diamond (x^{*}(t)\hsm -\hsm  y^{*}(t))},f(x^{*}(t),u)\rangle.\label{max-u}\end{eqnarray*}
\vspace{-.3cm}
\item[5.] Maximum condition on the upper level control.
\vspace{-.1cm}
\begin{eqnarray*}
&&q_H(t)+{\nu_L(t)\diamond  (x^*(t)- y^*(t))} \\
&& \hspace{.8cm}+ \textit{col} \Bigg(\sum_{j>1}  \nu_H^{{1j}}(t) \frac{y^{1*}(t)-y^{j*}(t)}{\|y^{1*}(t)-y^{j*}(t)\|} ;\vsm \\
&&\hspace{1.8cm} \ldots; \sum_{j<i}   \nu_H^{{ij}}(t) \frac{y^{i*}(t)-y^{j*}(t)}{\|y^{i*}(t)-y^{j*}(t)\|}\vsm \\
&& \hspace{2cm} + \sum_{j>i}   \nu_H^{{ij}}(t) \frac{y^{i*}(t)-y^{j*}(t)}{\|y^{i*}(t)-y^{j*}(t)\|};\vsm \\
&& \hspace{2cm}\ldots; \sum_{j<N}  \nu_H^{{Nj}}(t) \frac{y^{N*}(t)-y^{j*}(t)}{\|y^{N*}(t)-y^{j*}(t)\|} \Bigg)\vsm\vsm \\
&&\dis\hspace{1.cm} {\in \prod_{i=1}^N \alpha^i \partial^C_{v^i} \varphi^i(v^{i*}(t)) - \prod_{i=1}^N N_{{\cal V}^i}(v^{i*}(t)).}\label{max-vw}
\end{eqnarray*}
\vspace{-.1cm}

Here, each $\partial^C_{v^i} \varphi^i(v^{i*}(t))$, $i=1,\ldots,N$, is defined as below.
\end{itemize}

For $i=1,\ldots, N$, denote by $\Psi^i (y^i)$ the sets of optimal solutions to the lower level problems $(P_L^i)$, and by $H_L^i$ the Pontryagin-Hamiltonian functions given by
\begin{eqnarray*}
&&   H_L^i  \nonumber :=  \mu_L^i\langle x^i-y^i,v^i \rangle +\hspace{-.5cm}\sup\limits_{z \in - N^{M^i}_{D^i+y^i}(x^i)} \hspace{-.4cm} \{\langle p_L^i\hsm -\hsm \mu_L^i(x^i\hsm -\hsm y^i), z\rangle \}  \\
&& \hspace{1cm} \nonumber  +\sup_{u^i\in U^i}\{\langle p_L^i\hsm -\hsm \mu_L^i(x^i\hsm -\hsm y^i), f^i(x^i, u^i)\rangle\hsm -\hsm\bar\lambda^i \|u^i\|^2\} \\ && \hspace{.5cm}  + \left\langle p_H^i\hsm+\hsm \sum_{j<i} \mu_H^{{ij}}\frac{y^i-y^j}{\|y^i-y^j\|} \hsm + \hsm\sum_{j>i} \mu_H^{{ij}}\frac{y^i-y^j}{\|y^i-y^j\|} , v^i\right\rangle \nonumber
\end{eqnarray*}

\noindent for  $y^i,\, x^i,\,v^i,\, p_H^i,\, p_L^i $  in $\R^{2}$, $u^i\in \R^2$, $\mu_H^{ij},\mu_L^i, \bar \lambda^i\in \R_+$. Then, we have \begin{eqnarray*}\partial^C_{v^i} \varphi^i(v^{i*})& \hsm\hsm := \hsm\hsm & \textrm{co}\hsm\hsm\bigcup_{x^i\in\Psi^i (y^{i*})}\hsm\hsm \{\zeta^i \hsm\in\hsm L^2([0,T]\hsm :\R^{2})\hsm :\hsm  (\mathcal{A}^i) \textrm{ holds} \}. \end{eqnarray*}
The relation $(\mathcal{A}^i)$ is defined as follow (the dependence on $t$ is discarded for simplicity): there exists $(p_H^i,p_L^i, \mu_H^{ij},\mu_L^i,\bar \lambda^i) \in  AC \hsm \times \hsm AC  \hsm \times \hsm BV  \hsm \times \hsm BV  \hsm \times \hsm\R_+$ such that
\begin{itemize}
\item[a)] $\|(p_H^i,p_L^i)\|_{L^\infty}+ \|(\mu_H^{ij},\mu_L^i)\|_{TV} + \bar\lambda^i \neq 0$.\vspace{.1cm}
\item[b)] $\mu_H^{ij}$, $\mu_L^i $ are non-increasing, and constant, respectively, on $ \{t\hsm\in\hsm[0,T]\hsm:\hsm \|y^i-y^j\|\hsm >\hsm 2R, i\ne j, i,j\hsm =\hsm 1,\ldots,N\}$ and $\{t\hsm\in\hsm[0,T]\hsm:\hsm \|y^i- x^i\|\hsm <\hsm R\}$.\vspace{.1cm}
\item[c)] $ (\dot y^i,\dot x^i,-\dot{p}_H^i,-\dot{p}_L^i)\in \partial^C_{(p_H^i,p_L^i,y^i, x^i)}H_L^i\mbox{ a.e. } $\vspace{.2cm}
\item[d)] $\bar \lambda^i\zeta^i + p_H^i+\mu_L^i(x^i- y^i) $

$\qquad\quad  +\hsm \sum\limits_{j<i} \mu_H^{{ij}}\frac{y^i-y^j}{\|y^i-y^j\|}\hsm  +\hsm \sum\limits_{j>i} \mu_H^{{ij}}\frac{y^i-y^j}{\|y^i-y^j\|} \in - N_{\mathcal{V}^i}(v^i).	$\vspace{.2cm}
\item[e)]  $ p_H^i(0) \in \R^{2}$, \vspace{.1cm}

$ p_L^i(0)\hsm\in\hsm N_{D+y_0^i}(x^i(0))\hsm  +\hsm \mu_L^i(0) (x^i(0)\hsm -\hsm y_0^i)$,\vspace{.1cm}

$p_L^i(T) \hsm =\hsm  \mu_L^i(T) (x^i(T)\hsm -\hsm y^i(T))$,

\vspace{-.5cm}
 \begin{align*} \hspace{-.6cm} p_H^i(T)\hsm= & - \mu_L^i (T) (x^i(T) -y^i(T)) \\ &\hspace{1cm} -  \sum_{j<i}   \mu_H^{{ij}}(T) \frac{y^{i}(T)-y^{j}(T)}{\|y^{i}(T)-y^{j}(T)\|}\vsm \\ & \hspace{1.5cm}- \sum_{j>i}   \mu_H^{{ij}}(T) \frac{y^{i}(T)-y^{j}(T)}{\|y^{i}(T)-y^{j}(T)\|}. \end{align*}
\end{itemize}
\end{theorem}

\section{Brief outline of the key ideas of the proof}\label{sec:proof} The two main obstacles encountered when dealing with a bilevel control sweeping problem in the form of $(P_H)$-$(P_L)$ is first the discontinuity of the normal cone w.r.t. the state variable, and second the extra boundary constraint \eqref{low-level feasibility}. Indeed, the normal cone, being part of the dynamics of the lower level problem, lacks the Lipschitz property (with respect to the state variable) crucial to establish the standard necessary optimality condition. On the other hand, in the presence of condition \eqref{low-level feasibility}, the constraint qualifications, such as Mangasarian-Fromovitz or linear independence constraint qualifications, are likely to be violated and, thus, entailing a degeneracy phenomena.

In order to encounter these two challenges, we base the proof of Theorem \ref{theorem:main thm} mainly on two key pillars: the smooth approximation of the truncated normal cone by a sequence of Lipschitz functions in the state, and the flattening of the bilevel structure, under the partial calmness condition, by penalizing the cost of the upper level problem with the problematic term, i.e. the value function $\varphi$ representing condition {\eqref{low-level feasibility}}. The main ideas of the proof are based on the ones in \cite{KhalilPereiraIFAC2020}, however, adapted to our context, mainly in what concerns the consideration of $N\ge 2$ participating disks and the non-overlapping between them, while in \cite{KhalilPereiraIFAC2020}, only one disk is considered and therefore no overlapping is involved.

\section{Crowd motion example}\label{sec:examples}
\vspace{.1cm}
\begin{center}
\includegraphics[width=7cm, height=4.5cm]{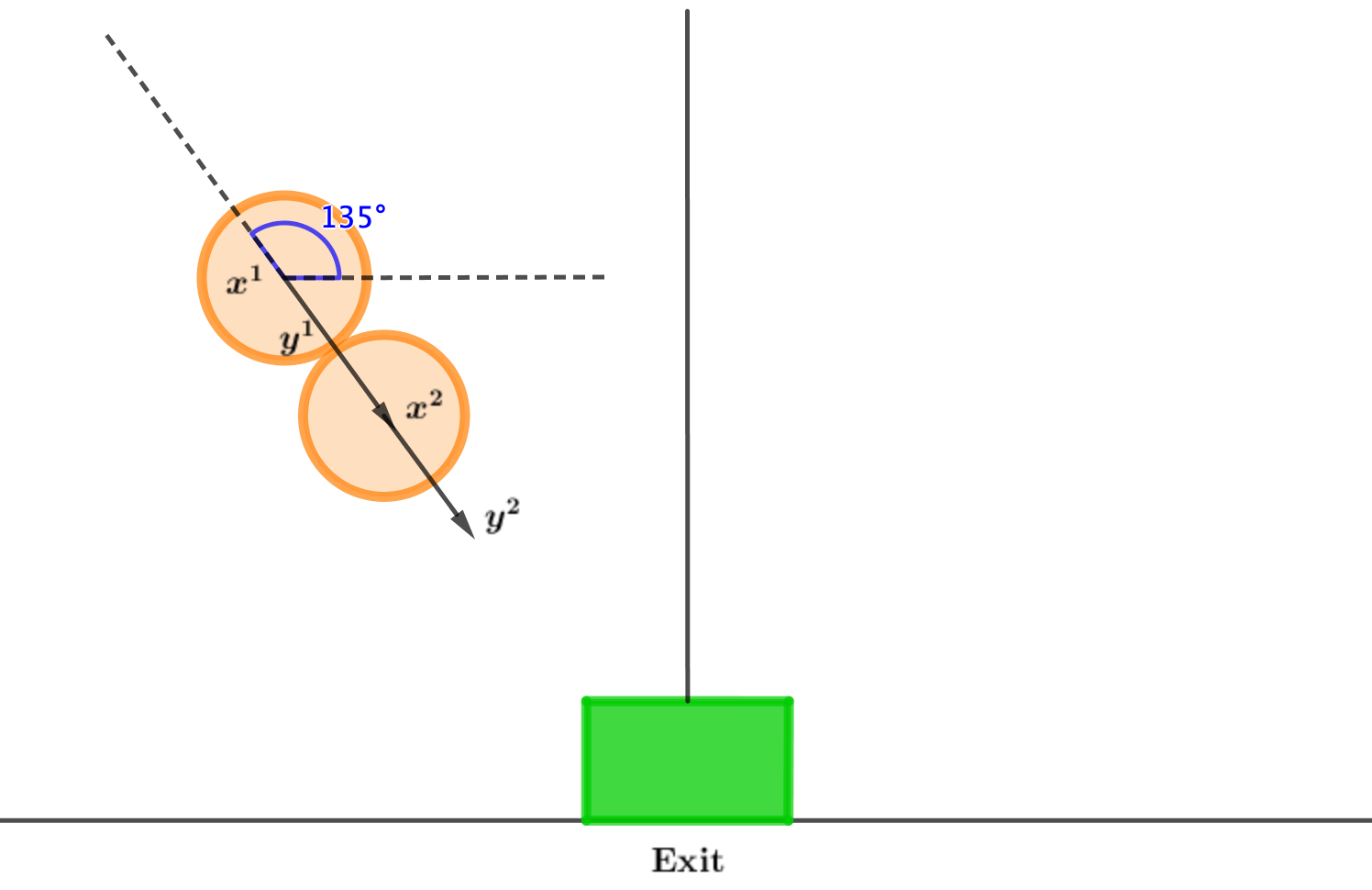}
\end{center}
In this example and to avoid confusion, we adopt the subscript notation when referring to the dynamics, controls, trajectories, and the state constraints, while keeping the superscript notation for the multipliers and the Hamiltonians.

We take two disks $D_i$, $i\hsm =\hsm 1,2$ of the same radius $R\hsm=\hsm3$. At the upper level problem, given a fixed final time $T\hsm =\hsm 6$, we consider the controls $v_1$, and $v_2$ taking values in {$V_i\hsm =\hsm \{ a (\cos(\theta_i),\sin(\theta_i))^T\hsm: \hsm a\hsm\in\hsm [-{\bf V},{\bf V}]\}$, with ${\bf V}\hsm =\hsm10 \sqrt{2}$}.  Vectors $y_1(t)$, and $y_2(t)$ defining the linear directions of disks, respectively, $D_1$, and $D_2$, satisfy the following dynamics
\begin{eqnarray*} &\dot{y}_1(t)= v_1=  -\bar v_1(\cos\theta_1,\sin\theta_1), \\ &  \; \dot{y}_2(t)=v_2=  -\bar v_2(\cos\theta_2,\sin\theta_2), \end{eqnarray*}
where $\theta_1$, and $\theta_2$ are the smallest positive angles formed by the positive $x$-axis with $y_1$, and $y_2$, respectively.

The disks keep the same direction until the final time, which entails that $\dot y_1(t)\hsm =\hsm \dot y_2(t)$ for all $t\in [0,6]$, and, thus, $\bar v_1\hsm =\hsm \bar v_2$.

At the upper level, we minimize the sum of the distances of the two disks to the exit at the origin:
\begin{equation*} \label{cost functio:high}\text{min } J_H:=\frac{1}{2} \left(  \|y_1(6)\|^2+  \|y_2(6)\|^2\right) .  \end{equation*}

Consider $y^0_1\hsm=\hsm(-48-3\sqrt{2},48\hsm +\hsm 3\sqrt{2})^T$, $ y^0_2\hsm =\hsm (-48,48^T)$, and $\theta_1\hsm =\hsm \theta_2\hsm =\hsm 135^\circ $. In this case, at the initial time $t\hsm =\hsm 0$, the two disks are in contact (indeed $\|y_1^0-y_2^0\|\hsm = \hsm 2R\hsm =\hsm 6$), and  $(cos\theta_i,\sin\theta_i)= {\bf v} $, for $i\hsm =\hsm 1,2$, where ${\bf v}= (-\frac{\sqrt{2}}{2},\frac{\sqrt{2}}{2})^T$.

At the lower level, we consider two scalar control functions $u_i(\cdot)$, $i\hsm =\hsm 1,2$ taking values in $U_i\hsm =\hsm [0,1]$ corresponding, to the efforts of the population of disks $D_i \hsm +\hsm y_i$. In each disk, the population is represented by its position $x_i$, $i\hsm =\hsm 1,2$. Denote by $t_a$ the contact time of $x_i$ with the boundary of the moving disks $D_i+y_i$. Moreover, once the population hits the boundary of the corresponding moving disk, it does not leave the boundary until the final time. The (average) dynamics of the populations are, for all $t\in [0,6]$:
\begin{eqnarray*}
	& \dot{x}_1(t)= - 8u_1(t)x_1(t) - \beta_1(t)u_0(x_1(t)-y_1(t)) \\ &  \dot{x}_2(t)= - 8u_2(t)x_2(t) - \beta_2(t)u_0(x_2(t)-y_2(t))
\end{eqnarray*}
being

$ \beta_1(t)\hsm=\hsm \begin{cases}
	\frac{M_1}{3} & \forall\; t\hsm\in\hsm [t_a,6]  \\ 0 &\forall \; t\hsm\in\hsm [0,t_a)
\end{cases}
$, $ \beta_2(t)\hsm=\hsm \begin{cases}
	\frac{M_2}{3} & \forall\; t\hsm\in\hsm [t_a,6]  \\ 0 &\forall\; t\hsm \in\hsm [0,t_a),
\end{cases}$
\vskip.5ex\noindent
$u_0\hsm\in\hsm [0,1]$ a control, and $M_1\hsm >\hsm 0$, and $M_2\hsm >\hsm 0$ are such that H5 is satisfied. Thus, we let $M_1\hsm =\hsm M_2\hsm =\hsm 6$.

The objective functions to minimize at the lower level are the control efforts of the populations to stay confined to their moving disks $D_1+y_1$ and $D_2+y_2$
\[\text{min }\int_{0}^{6} \|u_1(t)\|^2 dt, \quad \text{min }\int_{0}^{6} \|u_2(t)\|^2 dt.\]

Consider also for simplicity $(x_1^0,x_2^0)=(y_1^0,y_2^0)$, meaning that at the initial time $t=0$, each population is in the interior of its disk.

As follows from the problem formulation, the leader/follower solution concept entails the fact that the lower level problems relations constitute constraints to the upper level problem and there is no relation order on the lower level performance values.
This, together with the geometric insight dictate that $D_1$ and $ D_2 $ are driven from their initial position by the upper level dynamics along the line $(\alpha {\bf v}\hsm : \alpha\hsm\in\hsm\R \}$ to a final time position minimizing $J_H$ while enabling the lower level control problem constraints to be satisfied. On the other hand, the lower level controls should be such that $x^*_i$ stays as much as possible on the boundary of the disks $D_i+ y^*_i$, $i\hsm =\hsm 1,2$. It is straightforward to conclude that $v_1^*(t)\hsm =\hsm v_2^*(t)$ $\forall\, t\hsm\in\hsm[0,6]$ and that $\exists\, t_a,\,t_b $ such that $ 0\hsm < \hsm t_a\hsm<\hsm t_b\hsm<\hsm 6 $, and $\bar v^*\hsm >\hsm 0$ so that the optimal controls are
{\begin{eqnarray} \label{high-level-optimal control}\qquad\quad v_i^*(t)\hsm\hsm\hsm &= & \hsm\hsm\hsm
\left\{\begin{array}{ll} - \bar v^* {\bf v} & \mbox{if }t\hsm\in\hsm [0,t_b]\vspace{.1cm}\\
-(8\gamma_2(t)\hsm +\hsm 6){\bf v} & \mbox{if }t\hsm\in\hsm [t_b,6]\end{array}\right.\\
 \label{low-level-optimal control_1} \qquad\quad u_i^*(t)\hsm\hsm\hsm &= & \hsm\hsm\hsm
\left\{\begin{array}{ll} 0 & \mbox{if }t\hsm\in\hsm [0,t_a]\vspace{.1cm}\\
\dis \frac{\bar v^*-6}{8\|x_i^*(t)\|} & \mbox{if }t\hsm\in\hsm [t_a,t_b)\end{array}\right.\vspace{.1cm}\\
\label{low-level-optimal control_2}\qquad \forall\; t\hsm\in\hsm [t_b,6],\hspace{-.5cm}&&\hsm\hsm\hsm
\left\{\begin{array}{l} \dis u_1^*(t) =\frac{\|x_2^*(t)\|}{\|x_2^*(t)\|+6}\\
\dis u_2^*(t)= 1 \end{array}\right.
\end{eqnarray}}
where $\gamma_2$ is such that $x_2^*(t)\hsm =\hsm \gamma_2(t){\bf v}$ (note that $\|x^*_2\|\hsm =\hsm\gamma_2(t)$, whenever $\gamma_2\hsm\geq\hsm 0$). Observe that, for $i\hsm =\hsm 1,2$,
\begin{itemize}
\item The value of $\bar v^*$ on $[t_b,6]$ in (\ref{high-level-optimal control}) is required because $ u_2^*(t)\hsm\leq \hsm 1$, and $\dot x_2^*(t)\hsm =\hsm v_2^*(t) $, $\forall\, t\hsm\in\hsm [0,6]$.
\item In (\ref{low-level-optimal control_1}), the values of $u_i^*(t)$ on $[0,t_a]$ ensure that $x_i^*$ hits the boundary of $D_i+y_i^*$ as fast as possible, i.e., at $t\hsm =\hsm t_a$, which can be computed from $ \bar v^*t_a\hsm =\hsm 3$, and the ones in $(t_a, t_b]$, entails that $\dot x_i^*(t)=\dot y_i^*(t)$.
\item The values of $u_i^*(t)$ on $[t_b,6]$ in (\ref{low-level-optimal control_2}) reflect the upper bound on the values of $u_i^*$, and the relation $\dot x_1^*(t)\hsm =\hsm \dot x_2^*(t)$. Remark here that we always have $u_1^*(t)\hsm <\hsm 1$, $\forall\, t\hsm \in\hsm[0,6]$, and that $t_b$ and $\bar v^* $ can be related by the second control value expression in (\ref{low-level-optimal control_1}) for $i\hsm =\hsm 2$, i.e., $ \bar v^*\hsm =\hsm 8\|x_2^*(t_b)\|\hsm +\hsm 6$.
\end{itemize}
Now, we have that $\gamma_2(t)\hsm = \hsm 48\sqrt{2}\hsm +\hsm 3\hsm -\hsm \bar v^*t$, and, by solving $\dot\gamma_2(t)= -8\gamma_2-6$ on $[t_b,6]$, we obtain, for this interval, $ \gamma_2(t)\hsm =\hsm\frac{1}{8}(\bar v^* e^{-8(t-t_b)}- 6)$. From the continuity of $\gamma_2$ at $t_b$, we find that $\frac{1}{8}\bar v^*(8t_b+ 1)=48\sqrt{2}+ 3.75$. Moreover, by integrating $\bar v_2^* $, the value of $y_2(6)$ is obtained. After eliminating $\bar v^*$ with the previous equality we conclude that $$y_2(6)\hsm=\hsm \left(-3.75+\frac{48\sqrt{2}+3.750}{8t_b+1}e^{-8(6-t_b)}\right){\bf v}. $$ By noting that $y^*_1(t)\hsm =\hsm y^*_2(t)\hsm +\hsm 6{\bf v}$, and, by optimizing $J_H $ on $t_b$ as a function of $t_b$ alone, we obtain $t_b\hsm \approx \hsm 5.915 $, and, thus $\bar v^*\hsm\approx\hsm 11.860$, and $t_a\hsm \approx \hsm 0.253 $.

\medskip
Now, we show that the solution above satisfies the necessary optimality conditions of Theorem \ref{theorem:main thm}. We first notice that the assumptions stated in H1-H7 are satisfied. In what concerns checking the optimality conditions, we notice that there is a certain symmetry between the control processes $(y_1^*, v_1^*, x_1^*,u_1^*)$, and $(y_2^*, v_2^*, x_2^*,u_2^*)$ on the interval $[0,t_b]$. This symmetry is broken on $[t_b,6]$ due to the fact while $u_2^*$ is on the boundary of its constraints set that is not the case for $u_1^*$. Thus, the brief outline of the verification of the optimality conditions will be made with respect  to the former.

By expliciting the optimality conditions of Theorem \ref{theorem:main thm} at the computed $(v_1^*,v_2^*,u_1^*,u_2^*)$, we notice that the optimal control at the upper level $(v_1^*,v_2^*)$ is indeed in the interior of $V_1\times V_2$, otherwise $(v_1,v_2) \in \text{bd}(V_1\times V_2)$, equivalently, and either $v_1^*=v_2^*=0$ with the overall system at rest, or $v_1^*=v_2^*=10\sqrt{2}{\bf v}$. In either case the control process is obviously not optimal.

Then, the maximum condition 5. for the upper level problem yields
\vsm\vsm
\begin{eqnarray}  \nonumber
&&\hspace{-.3cm} q_H^2+\nu_L^2 (x_2^*- y_2^*) + \nu_H^{21} \frac{y^{*}_2-y_1^{*}}{\|y^{*}_2-y_1^{*}\|} \\
&&\hspace{.3cm}= -\frac{\alpha^2}{\bar\lambda^2} \bigg[ p_H^2\hsm+\hsm\mu_L^2(x_2^*-y_2^*) \hsm+\hsm\mu_H^{21} \frac{y^{*}_2-y_1^{*}}{\|y^{*}_2-y_1^{*}\|} \bigg] \label{eq:max cond lower level 1 - example}
\end{eqnarray}
\vsm 
which is a consequence of condition d) by remarking that
\begin{eqnarray*}  && \hspace{-.5cm}\partial_{v_2} \varphi_2 (v_2^*)\hsm\in\hsm \left\{-\frac{1}{\bar\lambda_2} \bigg[ p_H^2\hsm+\hsm\mu_L^2(x_2^*\hsm -\hsm y_2^*) \hsm+\hsm\mu_H^{21} \frac{y^{*}_2\hsm -\hsm y_1^{*}}{\|y^{*}_2\hsm -\hsm y_1^{*}\|} \bigg]\hsm \right\}.
\end{eqnarray*}
A similar analysis can be made with $q_H^1,\nu_L^1,\mu_L^1,\partial_{v_1} \varphi_1 (v_1^*)$.

In order to check whether the maximum condition is satisfied by the proposed optimal control $u_2^*$, we differentiate $H^2_H$ with respect to $u_2$, to obtain the expression 
\vsm
$$-8\bar x_2^*(\langle q_L^2,{\bf v}\rangle-\frac{1}{2}\nu_L^2 \|x_2^*- y_2^*\|)- 2\alpha^2 u_2^*,$$ 
\vsm
where $\bar x_2^*$ is such that $ x_2^*\hsm = \hsm\bar x_2^*{\bf v}$,
By concatenating the various segments and taking into account the system and adjoint dynamics, their boundary conditions, including the one articulating both levels, the state constraints, and their associated multipliers, we conclude that there exist multipliers for which $\nabla H^2_H\hsm \leq\hsm 0$ on $[0,t_a]$ whenever $ \|x_2^*\hsm -\hsm y_2^*\|\hsm <\hsm 3$, $\nabla H^2_H\hsm =\hsm 0$ on $[t_a,t_b]$ whenever $ \|x_2^*\hsm -\hsm y_2^*\|\hsm=\hsm 3$, being $H_H^2$ maximized by the proposed $u_2^*$, and $\nabla H^2_H\hsm \geq\hsm 0$ again when the low level state constraint is active, forcing $u_2^*\hsm=\hsm1$. Similar reasoning can be made for $u_1^*$.

Now, by evaluating \eqref{eq:max cond lower level 1 - example} at the final time $T=6$, and assuming (without loss of generality) that all the measures at $T=6$ are equal to zero, we obtain, from the boundary conditions (condition 3.) and from $y^*_1(6)\neq 0$ (also $y^*_2(6)\neq 0$) that $\lambda\hsm=\hsm0$ (yielding that $\alpha^1\hsm=\hsm\alpha^2\hsm=\hsm0$). By substituting again in \eqref{eq:max cond lower level 1 - example}, but now for a.e. $t\hsm\in\hsm  [t_a,6]$ and by taking the norm, we obtain that $\|(q_H^1,q_H^2)\| \hsm =\hsm \|(\nu_L^1,\nu_L^2)\| \hsm =\hsm\|\nu_H^{12}\|=0$. By replacing in the maximum condition on the lower level problem (condition 4.), it is clear that $q_L^1$ and $q_L^2$ cannot be zero on $[t_a,6]$, otherwise, the maximum would be attained at $u_1^*\hsm =\hsm u_2^*=0$, contradicting the expressions of the optimal controls in \eqref{low-level-optimal control_1}, and \eqref{low-level-optimal control_2} on $[0.254,6]$. This confirms the nontriviality of multipliers. Hence, the necessary optimality conditions postulated in Theorem \ref{theorem:main thm} are satisfied at $(\bar v_1^*,\bar v_2^*,u_1^*,u_2^*)$.

\section{Conclusion}\label{sec:Conclusions}
In this work, we investigated the structured crowd motion model formulated as a bilevel control sweeping problem. We established the necessary optimality conditions when several participants are involved. To better understand those conditions, we illustrate an example with only two participants, for which we find the optimal solution explicitly. We showed that the necessary optimality conditions are satisfied for the computed optimal solution.

For the future research, we are interested in investigating a setting of problem where a sweeping process is incorporated in the dynamics of the upper level as well. This occurs for instance if we consider $m$ groups of crowd motion models, where each group $i\in\{1,\ldots,m\}$ is following its own agent, and structured into $N^i$ subgroups. At the upper level, while each agent tends to reach its target with the minimal possible time, the sweeping process arises so that each agent guarantees the non-overlapping of his/her own group with the others. At the lower level, the subgroups will follow their own agents with the minimum effort while keeping a safe distance with the agent and with the other subgroups of the same group. This translates into a sweeping process in the corresponding dynamics at the lower level too.


\section*{Acknowledgment}
Tan  H.  Cao  acknowledges the  support  of  the National Research Foundation of Korea grant funded by the Korea Government (MIST) NRF-2020R1F1A1A01071015.

N.T. Khalil and F.L. Pereira acknowledge the support of FCT R\&D Unit SYSTEC - POCI-01-0145-FEDER-006933, and MAGIC - POCI-01-0145- FEDER-032485 - funded by ERDF - COMPETE2020 - FCT/MEC - PT2020, STRIDE - NORTE-01-0145-FEDER-000033.

Boris S. Mordukhovich acknowledges the sup-port  of  the  USA  National  Science  Foundation  under  grants DMS-1007132  and  DMS-1512846,  by  the  USA  Air  Force Office  of  Scientific  Research  grant  $\#$15RT0462,  and  by  the Australian  Research  Council  under  Discovery  Project  DP-190100555.

Dao Nguyen acknowledges the support of the USA National Science Foundation under grant DMS-1512846,  and  by  the  USA  Air  Force  Office  of Scientific Research grant $\#$15RT0462.


\bibliographystyle{apalike}

\end{document}